\documentclass[leqno]{article}

\usepackage{amsmath}
\usepackage{amsthm}
\usepackage{amsfonts}
\begin {document}
\topmargin= -.2in \baselineskip=20pt

\newtheorem{theorem}{Theorem}[section]
\newtheorem{proposition}[theorem]{Proposition}
\newtheorem{lemma}[theorem]{Lemma}

\theoremstyle{remark}
\newtheorem{remark}[theorem]{Remark}
\def\mathscr{\mathcal}

\title {A Class of Incomplete Character Sums
\thanks{We would like to thank the referee for suggestion of applying tensor induction
to higher rank sheaves. Originally we only treat rank one sheaves
using transfer. We would like to thank Jos\'e Alves Oliveira and Lucas Reis for pointing out an error in 
the old version of Theorem 4.6. The research of Lei Fu is supported by NSFC.}}

\author {Lei Fu\\
{\small Chern Institute of Mathematics and LPMC, Nankai University,
Tianjin 300071, P. R. China}\\
{\small leifu@nankai.edu.cn}\\{}\\
Daqing Wan\\
{\small Department of Mathematics, University of California, Irvine,
CA
92697, USA}\\
{\small dwan@math.uci.edu}}
\date{}

\maketitle

\begin{abstract}
Using $\ell$-adic cohomology of tensor inductions of lisse
$\overline{\mathbb Q}_\ell$-sheaves, we study a class of incomplete
character sums.



\end{abstract}

\section*{Introduction}

Through out this paper, $p$ is a fixed prime number, $\mathbb F_p$
is a finite field with $p$ elements, and $\mathbb F$ is an algebraic
closure of $\mathbb F_p$. For any power $q$ of $p$, let $\mathbb
F_q$ be the subfield of $\mathbb F$ with $q$ elements. Let $\ell$ be
a prime number distinct from $p$, let $\chi_1,\ldots, \chi_k:\mathbb
F_{q^d}^\ast\to\overline{\mathbb Q}_\ell^\ast$ be a family of
multiplicative characters  and $\psi:\mathbb
F_{q^d}\to\overline{\mathbb Q}_\ell^\ast$ an additive character on
an extension field $\mathbb F_{q^d}$ of $\mathbb F_q$. We extend
$\chi_i$ $(i=1,\ldots, k)$ to $\mathbb F_{q^d}$ by setting
$\chi_i(0)=0$. Let $f_1(t_1,\ldots, t_n),\ldots,
f_k(t_1,\ldots,t_n), f_{k+1}(t_1,\ldots, t_n)\in \mathbb
F_{q^d}[t_1,\ldots, t_n]$ be a family of polynomials with
coefficients in $\mathbb F_{q^d}$. Motivated by a number of
applications in \cite{W}, we are interested in  estimating the
following type of incomplete character sum
$$S'=\sum_{x_1,\ldots, x_n\in\mathbb F_q}
\chi_1\big(f_1(x_1,\ldots, x_n)\big)\cdots
\chi_k\big(f_k(x_1,\ldots, x_n)\big)\psi\big(f_{k+1}(x_1,\ldots,
x_n)\big),$$ where the summation is over those $x_1,\ldots, x_n$ in
the subfield $\mathbb F_q$ of $\mathbb F_{q^d}$. Note that for the
classical complete character sum
$$S=\sum_{x_1,\ldots, x_n\in\mathbb F_{q^d}}
\chi_1\big(f_1(x_1,\ldots, x_n)\big)\cdots
\chi_k\big(f_k(x_1,\ldots, x_n)\big)\psi\big(f_{k+1}(x_1,\ldots,
x_n)\big),$$ the summation is over all $x_1,\ldots, x_n$ in the
field $\mathbb F_{q^d}$. Under suitable hypothesis, one can give a
sharp estimate of the following form for the complete sum
$$|S| \leq C(n,\mathrm{deg}(f_i)) {q}^{\frac{nd}{2}},$$
where $C(n,\mathrm{deg}(f_i))$ is a constant depending only on $n$
and the degrees of the functions $f_i$'s. We believe that there is
also a sharp estimate for the incomplete sum $S'$. Namely, under
suitable hypothesis, there should be an estimate of the form
$$|S'| \leq C(n, d,\mathrm{deg}(f_i)) {q}^{\frac{n}{2}},$$
where $C(n, d,\mathrm{deg}(f_i))$ is a constant depending only on
$n, d$ and the degrees of $f_i$. This is indeed true in the one
variable case $n=1$, see \cite{W}. In the present paper, we use
tensor induction and $\ell$-adic cohomology theory to study this
problem.

Let us recall the general method of studying (complete) character
sum using $\ell$-adic cohomology theory. Let $X_0$ be a separated
scheme of finite type over $\mathbb F_q$, let $X=X_0\otimes_{\mathbb
F_q}\mathbb F$ be the base change of $X_0$, let $\mathscr L_0$ be a
lisse $\overline{\mathbb Q}_\ell$-sheaf on $X_0$, and let $\mathscr
L$ be the inverse image of $\mathscr L_0$ in $X$. The family of
complete character sums corresponding to $\mathscr L_0$ is the
family of sums
$$S_m(X_0,\mathscr L_0)=\sum_{x\in X_0(\mathbb F_{q^m})}
\mathrm{Tr}(F_x,\mathscr L_{0,\bar x}),$$ where $X_0(\mathbb
F_{q^m})=\mathrm{Hom}_{\mathrm{Spec}\,\mathbb
F_q}(\mathrm{Spec}\,\mathbb F_{q^m},X_0)$ is the set of $\mathbb
F_{q^m}$-points in $X_0$, $\mathscr L_{0,\bar x}$ is the stalk of
$\mathscr L_0$ at a geometric point over $x$, and
$F_x\in\mathrm{Gal}(\mathbb F/\mathbb F_{q^m})$ is the geometric
Frobenius element at $x$, that is, the inverse of the Frobenius
substitution. By Grothendieck's trace formula (\cite[Rapport
3.2]{SGA41/2}), we have
$$S_m(X_0,\mathscr L_0)=\sum_{i=0}^{2\mathrm{dim}(X)}(-1)^i \mathrm{Tr}\big(F^m,
H^i_c(X,\mathscr L)\big),$$ where $H^i_c(X,\mathscr L)$ are
$\ell$-adic cohomology groups with compact support and $F$ is the
geometric Frobenius correspondence. Deligne's theorem
\cite[Corollarie 3.3.4]{D1} gives an estimate for the Archimedean
absolute values of the eigenvalues of $F$ on $H^i_c(X,\mathscr L)$.
If we have some information about the Betti numbers
$\mathrm{dim}\,H^i_c(X,\mathscr L)$, then we can get an estimate for
$S_m(X_0,\mathscr L_0)$.

Let $\mathbb F_{q^d}$ be an extension of $\mathbb F_q$ of degree
$d$, let $X_1=X_0\otimes_{\mathbb F_q}\mathbb F_{q^d}$ be the base
change of $X_0$, and let $\mathscr L_1$ be a lisse
$\overline{\mathbb Q}_\ell$-sheaf on $X_1$. For any $\mathbb
F_q$-point $x:\mathrm{Spec}\,\mathbb F_q\to X_0$ of $X_0$, let
$x_1:\mathrm{Spec}\,\mathbb F_{q^d}\to X_1$ be its base change. It
is a $\mathbb F_{q^d}$-point of $X_1$. In this paper, we are
concerned with the following incomplete character sum
$$S_1'(X_1, \mathscr L_1)=\sum_{x\in X_0(\mathbb
F_q)}\mathrm{Tr}(F_{x_1}, \mathscr L_{1,\bar x_1}),$$ where the
summation is over those $\mathbb F_{q^d}$-points of $X_1$ arising
from $\mathbb F_q$-points of $X_0$ by base change. In section 1, we
describe the tensor induction $\otimes\hbox{-{\rm Ind}}(\mathscr
L_1)$ of $\mathscr L_1$ (\cite[10.3-5]{K}), which is a lisse
$\overline{\mathbb Q}_\ell$-sheaf of rank $1$ on $X_0$. We show that
the above incomplete sum $S_1'(X_1,\mathscr L_1)$ coincides with the
following complete character sum
$$S_1(X_0,\otimes\hbox{-{\rm Ind}}(\mathscr L_1))=\sum_{x\in X_0(\mathbb
F_q)}\mathrm{Tr}(F_x, (\otimes\hbox{-{\rm Ind}}(\mathscr L_1))_{\bar
x}).$$ In principle, this reduces the study of our incomplete sum to
the study of the complete sum for the new sheaf $\otimes\hbox{-{\rm
Ind}}\mathscr (L_1)$. We describe explicitly the tensor induction
for Kummer type sheaves in \S 2, and for Artin-Schreier type sheaves
in \S 3. In general, it is hard to estimate the Betti numbers
$\mathrm{dim}\,H^i_c(X, \otimes\hbox{-{\rm Ind}}(\mathscr L_1))$. In
\S 4, we study the case where $\mathrm{dim}(X)=1$. In this case, the
Grothendieck-Ogg-Shararevich formula for Euler characteristic
numbers of sheaves on curves gives us more information about the
Betti numbers. We can thus get relatively complete results (Theorems
4.5-6).

\section{Tensor Induction}

Let $G$ be a pro-finite group, let $H$ be an open subgroup of $G$ of
finite index $d$, and let $\rho: H\to \mathrm{GL}(V)$ be a
(continuous) representation of $H$, where $V$ is a finitely
dimensional vector space over $\overline{\mathbb Q}_\ell$. Choose
representatives ${Hg_1,\ldots, Hg_d}$ for the family of right cosets
$H\backslash G$. For any $g\in G$, the set $\{Hg_1g,\ldots, Hg_dg\}$
is also a family of representatives of right cosets. So we have
$$Hg_ig=Hg_{\tau(i)}$$ for a permutation $i\mapsto \tau(i)$
of $\{1,\ldots, d\}$. Define $$\otimes\hbox{-{\rm Ind}}(\rho):G\to
\mathrm{GL}(\otimes^d V)$$ by
$$\Big((\otimes\hbox{-{\rm Ind}}(\rho)(g)\Big)(v_1\otimes \cdots \otimes
v_d)=\rho(g_1 g g^{-1}_{\tau(1)})(v_{\tau(1)})\otimes \cdots\otimes
\rho(g_d g g^{-1}_{\tau(d)})(v_{\tau(d)})$$ for any $v_1,\ldots,
v_d\in V$. One can show that $\otimes\hbox{-{\rm Ind}}(\rho)$ is a
(continuous) representation of $G$, and its isomorphism class is
independent of the choice of representatives of right cosets. We
call $\otimes\hbox{-{\rm Ind}}(\rho)$ the \emph{tensor induction} of
$\rho$. For more details about the tensor induction, see \cite[\S
13]{CR} and \cite[10.3-5]{K}. In the following, we summarize the
properties of tensor induction which are used in the paper.

\begin{proposition} Let $G$ be a pro-finite group, let $H$ be an open subgroup
of $G$ of finite index $d$, and let $\rho:H\to\mathrm{GL}(V)$ be a
representation of $H$.

(i) Suppose that $H$ is normal in $G$. Let $\{Hg_1,\ldots, Hg_d\}$
be a family of representatives for $H\backslash G$, and let
$\rho^{(i)}:H\to \mathrm{GL}(V)$ $(i=1,\ldots, d)$ be the
representations defined by $\rho^{(i)}(h)=\rho(g_i hg_i^{-1})$ for
any $h\in H$. Then for any $h\in H$, we have
$$\otimes\hbox{-{\rm Ind}}(\rho)(h)=\otimes_{i=1}^d \rho^{(i)}(h).$$

(ii) Suppose that $H$ is normal in $G$ and $G/H$ is a cyclic group.
For any $\sigma\in G$ such that $\sigma H$ is a generator of $G/H$,
we have
$$\mathrm{Tr}(\otimes\hbox{-{\rm Ind}}(\rho)(\sigma))= \mathrm{Tr}(\rho(\sigma^d)).$$
\end{proposition}

\begin{proof} (i) follows directly from the definition of the tensor
induction. To prove (ii), use $He, H\sigma, \cdots, H\sigma^{d-1}$
as representatives for right cosets. Let $\{e_1, \ldots, e_n\}$ be a
basis of $V$. Then $\{e_{j_1}\otimes \cdots\otimes e_{j_d}| 1\leq
j_1,\ldots, j_d\leq n\}$ is a basis of $\otimes^d V$. We have
$$\otimes\hbox{-{\rm Ind}}(\rho)(\sigma)(e_{i_1}\otimes \cdots\otimes
e_{i_d})=e_{i_2}\otimes\cdots\otimes e_{i_{d}}\otimes
\rho(\sigma^d)(e_{i_1}).$$ Suppose $\rho(\sigma^d)(e_i)=\sum_{j=1}^n
a_{ij}e_j$. In the expansion of $\otimes\hbox{-{\rm
Ind}}(\rho)(\sigma)(e_{i_1}\otimes \cdots\otimes e_{i_d})$ as a
linear combination of $\{e_{j_1}\otimes \cdots\otimes e_{j_d}| 1\leq
j_1,\ldots, j_d\leq n\}$, the coefficient of $e_{i_1}\otimes
\cdots\otimes e_{i_d}$ is nonzero only when $i_1, \ldots, i_d$ are
equal to a common $i$, and in this case, the coefficient is
$a_{ii}$. It follows that
$$\mathrm{Tr}(\otimes\hbox{-{\rm Ind}}(\rho)(\sigma))= \sum_i
a_{ii}= \mathrm{Tr}(\rho(\sigma^d)).$$
\end{proof}

Let $X_0$ be a scheme of finite type over $\mathbb F_q$, and let
$X_1=X_0\otimes_{\mathbb F_q}\mathbb F_{q^d}$ be its base change to
an extension field $\mathbb F_{q^d}$ of degree $d$ over $\mathbb
F_q$. Suppose $X_0$ is geometrically connected, that is, the base
change $X=X_0\otimes_{\mathbb F_q}\mathbb F$ is connected. Fix a
geometric point for $X$. Take its images in $X_0$ and $X_1$ as the
base points of $X_0$ and $X_1$, respectively, and let $\pi_1(X_0)$
and $\pi_1(X_1)$ be the \'etale fundamental groups with respect to
these base points. Then $\pi_1(X_1)$ is a normal subgroup of
$\pi_1(X_0)$. We have an isomorphism
$$\pi_1(X_0)/\pi_1(X_1)\stackrel \cong \to \mathrm{Gal}(\mathbb
F_{q^d}/\mathbb F_q),$$ and $\mathrm{Gal}(\mathbb F_{q^d}/\mathbb
F_q)$ is a cyclic group. Let $x:\mathrm{Spec}\,\mathbb F_{q^m}\to
X_0$ be an $\mathbb F_{q^m}$-points in $X_0$. We can talk about the
geometric Frobenius element $F_x$ in $\pi_1(X_0)$ corresponding to
$x$. It is defined up to conjugation in $\pi_1(X_0)$. Consider the
case where $m=1$. Then the image of $F_x$ in the quotient group
$\pi_1(X_0)/\pi_1(X_1)$ is a generator. Let $x_1:
\mathrm{Spec}\,\mathbb F_{q^d}\to X_1$ be the base change of $x$.
One can verify that $F_x^d$ and $F_{x_1}$ define the same conjugacy
class in $\pi_1(X_1)$.

Let $\mathscr L_1$ be a lisse $\overline{\mathbb Q}_\ell$-sheaf on
$X_1$ of rank $r$. It defines a representation
$$\rho_1: \pi_1(X_1)\to \mathrm{GL}(r, \overline{\mathbb Q}_\ell).$$
Regard $\pi_1(X_1)$ as a subgroup $\pi_1(X_0)$, and consider the
tensor induction
$$\otimes\hbox{-{\rm Ind}}(\rho_1): \pi_1(X_0)\to
\mathrm{GL}(r^d,\overline{\mathbb Q}_\ell).$$ We define the
\emph{tensor induction} $\otimes\hbox{-{\rm Ind}}(\mathscr L_1)$ of
$\mathscr L_1$ to be the lisse $\overline{\mathbb Q}_\ell$-sheaf on
$X_0$ of rank $r^d$ corresponding to the last representation. By the
above discussion and Proposition 1.1 (ii), if
$x:\mathrm{Spec}\,\mathbb F_q\to X_0$ is an $\mathbb F_q$-point in
$X_0$ and $x_1: \mathrm{Spec}\,\mathbb F_{q^d}\to X_1$ is its base
change, then we have
\begin{eqnarray}
\mathrm{Tr}\big(F_x, (\otimes\hbox{-{\rm Ind}}(\mathscr L_1))_{\bar
x}\big)=\mathrm{Tr}(F_{x_1},\mathscr L_{1,\bar x_1}).
\end{eqnarray}

\begin{proposition} Keep the notation above.

(i) The incomplete character sum
$$S_1'(X_1, \mathscr L_1)=\sum_{x\in X_0(\mathbb
F_q)}\mathrm{Tr}(F_{x_1}, \mathscr L_{1,\bar x_1})$$ coincides with
the complete character sum
$$S_1(X_0,\otimes\hbox{-{\rm Ind}}(\mathscr L_1))=\sum_{x\in X_0(\mathbb
F_q)}\mathrm{Tr}\big(F_x, (\otimes\hbox{-{\rm Ind}}(\mathscr
L_1))_{\bar x}\big).$$

(ii) Let $\rho_1:\pi_1(X_1)\to \mathrm{GL}(r,\overline{\mathbb
Q}_\ell)$ be the representation corresponding to the
$\overline{\mathbb Q}_\ell$-sheaf $\mathscr L_1$, and let $\sigma\in
\pi_1(X_0)$ be an element such that its image under the homomorphism
$\pi_1(X_0)\to\mathrm{Gal}(\mathbb F_{q^d}/\mathbb F_q)$ is the
Frobenius substitution. For each $0\leq i\leq d-1$, let
$\rho_1^{(i)}:\pi_1(X_1)\to \mathrm{GL}(r, \overline{\mathbb
Q}_\ell)$ be the representation defined by
$$\rho_1^{(i)}(g)=\rho_1(\sigma^{i}g\sigma^{-i}),$$ and let $\pi:X_1\to
X_0$ be the projection. Then the representation corresponding the
sheaf $\pi^\ast(\otimes\hbox{-{\rm Ind}}(\mathscr L_1))$ on $X_1$ is
$\otimes_{i=0}^{d-1} \rho_1^{(i)}$.

(iii) Let $\mathscr L_1$ and $\mathscr M_1$ be two lisse
$\overline{\mathbb Q}_\ell$-sheaves. Then we have
$$\big(\otimes\hbox{-{\rm Ind}}\big)(\mathscr L_1\otimes \mathscr M_1) \cong
\Big(\otimes\hbox{-{\rm Ind}}(\mathscr L_1)\Big)\otimes
\Big(\otimes\hbox{-{\rm Ind}}(\mathscr M_1)\Big).$$
\end{proposition}

\begin{proof}
(i) follows from Equation (1). (ii) follows from Proposition 1.1
(i). (iii) follows from the definition of tensor induction.
\end{proof}

\section{Tensor Induction of Kummer type sheaves}

The Kummer covering
$$[q-1]:\mathbb G_{m,\mathbb F_q}\to \mathbb G_{m,\mathbb F_q},
\quad x\mapsto x^{q-1}$$ on $\mathbb G_{m,\mathbb F_q}=\mathrm
{Spec}\,\mathbb F_q[t,t^{-1}]$ defines a $\mathbb F_q^\ast$-torsor
$$1\to \mathbb F_q^\ast\to \mathbb G_{m,\mathbb F_q}
\stackrel {[q-1]}\to \mathbb G_{m,\mathbb F_q}\to 1.$$ Let $\chi:
\mathbb F_q^\ast\to \overline {\mathbb Q}_\ell^\ast$ be a
multiplicative character. Pushing-forward the above torsor by
$\chi^{-1}$, we get a lisse $\overline {\mathbb Q}_\ell$-sheaf
${\mathscr K}_\chi$ on $\mathbb G_{m,\mathbb F_q}$ of rank 1. We
call ${\mathscr K}_\chi$ the \emph{Kummer sheaf} associated to
$\chi$. For any $x\in\mathbb G_{m,\mathbb F_q}(\mathbb
F_{q^m})=\mathbb F_{q^m}^\ast$, we have
$$\mathrm{Tr}(F_x,\mathscr K_{\chi,\bar x})=\chi(\mathrm{N}_{\mathbb F_{q^m}/\mathbb
F_q}(x)),$$ where $\mathrm{N}_{\mathbb F_{q^m}/\mathbb F_q}$ denotes
the norm for the field extension $\mathbb F_{q^m}/\mathbb F_q$.

Let $f\in\Gamma(X_0,\mathcal O_{X_0}^\ast)$ be a section of the
subsheaf of units in the structure sheaf $\mathcal O_{X_0}$. Then
the $\mathbb F_q$-algebra homomorphism
$$\mathbb F_q[t,t^{-1}]\to \Gamma(X_0,\mathcal O_{X_0}),\quad t\mapsto f$$
defines an $\mathbb F_q$-morphism $X_0\to\mathbb G_{m,\mathbb F_q}$
which we still denote by $f$. Denote $f^\ast \mathscr K_\chi$ by
$\mathscr K_{\chi,f}$. It is a lisse $\overline{\mathbb
Q}_\ell$-sheaf of rank one on $X_0$. For any $x\in X_0(\mathbb
F_{q^m})$, we have
$$\mathrm{Tr}(F_x,(\mathscr K_{\chi,f})_{\bar x})=\chi\big(\mathrm{N}_{\mathbb F_{q^m}/\mathbb
F_q}(f(x))\big).$$ For any $f_1,f_2, f\in\Gamma(X_0,\mathcal
O_{X_0}^\ast)$ and any multiplicative characters $\chi,
\chi_1,\chi_2: \mathbb F_q^\ast\to\overline{\mathbb Q}_\ell^\ast$,
we have
$$\mathscr K_{\chi,f_1f_2}\cong \mathscr K_{\chi,f_1}\otimes
\mathscr K_{\chi,f_2},\quad \mathscr K_{\chi_1\chi_2, f}\cong
\mathscr K_{\chi_1,f}\otimes \mathscr K_{\chi_2,f}.$$ Confer
\cite[Sommes trig. 1.4-8]{SGA41/2}.

Suppose that $X_0$ is a normal scheme. Then $\mathscr K_{\chi,f}$
can be described as follows. Let $K$ be the function field of $X_0$.
Fix a separable closure of $\overline {K}$ of $K(X_0)$. Let the base
point of $X_0$ be the canonical morphism $\mathrm{Spec}\,\overline
K\to X_0$. Then we have a canonical isomorphism
$$\pi_1(X_0)\cong \varprojlim_{L}\mathrm{Gal}(L/K),$$ where $L$ goes over those
finite Galois extensions of $K$ contained in $\overline K$ such that
the normalization of $X_0$ in $L$ is etale over $X_0$. Let $\xi\in
\overline K$ be a root of the polynomial $T^{q-1}-f.$ Then any root
of $T^{q-1}-f$ is of the form $a\xi$ for some $a\in\mathbb
F_q^\ast$,  $K[\xi]$ is a Galois extension, and the normalization of
$X_0$ in $K[\xi]$ is etale over $X_0$. It follows that we have a
canonical epimorphism
$$\pi_1(X_0)\to \mathrm{Gal}(K[\xi]/K).$$ On the other hand, we have
a canonical monomorphism
$$\mathrm{Gal}(K[\xi]/K)\hookrightarrow \mathbb F_q^\ast,\quad
g\mapsto \frac{g(\xi)}{\xi}.$$ The representation corresponding to
$\mathscr K_{\chi,f}$ is the composite
$$\pi_1(X_0)\to \mathrm{Gal}(K[\xi]/K) \hookrightarrow \mathbb F_q^\ast
\stackrel{\chi^{-1}}\to \overline {\mathbb Q}_\ell^\ast.$$

Let $\mathbb F_{q^d}$ be an extension of $\mathbb F_q$ of degree
$d$, let $\chi:\mathbb F_{q^d}^\ast\to \overline{\mathbb
Q}_\ell^\ast$ be a multiplicative character, let $f\in\mathbb
F_{q^d}(t_1,\ldots, t_n)$ be a rational function, and let $X_0$ be
an open subscheme of $\mathbb A^n_{\mathbb F_q}$ so that its inverse
image in $\mathbb A^n_{\mathbb F_{q^d}}$ is contained in the
complement of the union of hypersurfaces $\{f=0\}\cup\{f=\infty\}$.
Let $X_1$ be the base change of $X_0$ from $\mathbb F_q$ to $\mathbb
F_{q^d}$, and let $\pi:X_1\to X_0$ be the projection. Consider the
sheaf $\mathscr K_{\chi,f}$ on $X_1$. Let us describe the character
$\pi_1(X_1)\to\overline{\mathbb Q}_\ell^\ast$ corresponding to the
rank one lisse sheaf $\pi^\ast(\otimes\hbox{-{\rm Ind}}(\mathscr
K_{\chi,f}))$ on $X_1$. The function field of $X_0$ (resp. $X_1$) is
the rational function field $K_0=\mathbb F_q(t_1,\ldots, t_n)$
(resp. $K_1=\mathbb F_{q^d}(t_1,\ldots, t_n)$). Let $\xi$ be a root
of the polynomial $T^{q^d-1}-f$ in a separable closure $\overline
K_1$ of $K_1={\mathbb F_{q^d}(t_1,\ldots, t_n)}$. Then $\mathscr
K_{\chi,f}$ corresponds to the character
$$\pi_1(X_1)\to
\mathrm{Gal}(K_1[\xi]/K_1)\hookrightarrow \mathbb
F_{q^d}^\ast\stackrel{\chi^{-1}}\to \overline{\mathbb
Q}_\ell^\ast.$$ Denote this character by
$\rho_1:\pi_1(X_1)\to\overline{\mathbb Q}_\ell^\ast$. More
explicitly, we have
$$\rho_1(g)=\chi^{-1}\Big(\frac{g(\xi)}{\xi}\Big)$$ for any $g\in
\pi_1(X_1)$. Note that $\pi:X_1\to X_0$ is a Galois etale covering
space, and we have canonical isomorphisms
$$\mathrm{Aut}(X_1/X_0)^\circ\cong\mathrm{Gal}\Big(\mathbb F_{q^d}(t_1,\ldots,
t_n)/\mathbb F_{q}(t_1,\ldots, t_n)\Big)\cong \mathrm{Gal}(\mathbb
F_{q^d}/\mathbb F_q).$$ Choose an element $\sigma\in \pi_1(X_0)$ so
that its image in $\mathrm{Gal}\Big(\mathbb F_{q^d}(t_1,\ldots,
t_n)/\mathbb F_{q}(t_1,\ldots, t_n)\Big)$ is given by
$$\sigma(a)=a^q \hbox { for any } a\in\mathbb F_{q^d},\quad
\sigma(t_i)=t_i\; (i=1, \ldots, n).$$ Then the image of $\sigma$ in
$\pi_1(X_0)/\pi_1(X_1)$ is a generator. For any $0\leq i\leq d-1$,
let $\rho_1^{(i)}:\pi_1(X_1)\to\overline{\mathbb Q}_\ell^\ast$ be
the character defined by
$$\rho_1^{(i)}(g)=\rho_1(\sigma^ig\sigma^{-i})$$ for any $g\in
\pi_1(X_1)$. We have
\begin{eqnarray*}
\rho_1^{(i)}(g)&=&\rho_1(\sigma^{i}g\sigma^{-i})\\
&=& \chi^{-1}\Big(\frac{\sigma^{i}g\sigma^{-i}(\xi)}{\xi}\Big)\\
&=&\chi^{-1}\Big(\sigma^{i}\Big(\frac{g\sigma^{-i}(\xi)}{\sigma^{-i}(\xi)}\Big)\Big).
\end{eqnarray*}
Note that $\sigma^{-i}(\xi)$ and $g\sigma^{-i}(\xi)$ are roots of
the polynomial $T^{q^d-1}-\sigma^{-i}(f)$, and hence
$\frac{g\sigma^{-i}(\xi)}{\sigma^{-i}(\xi)}$ lies in $\mathbb
F_{q^d}^\ast$. We thus have
$$\sigma^{i}\Big(\frac{g\sigma^{-i}(\xi)}{\sigma^{-i}(\xi)}\Big)=
\Big(\frac{g\sigma^{-i}(\xi)}{\sigma^{-i}(\xi)}\Big)^{q^i}.$$ So we
have
$$\rho_1^{(i)}(g)=(\chi^{q^i})^{-1}\Big(\frac{g\sigma^{-i}(\xi)}{\sigma^{-i}(\xi)}\Big)$$
for any $g\in \pi_1(X_1)$. This shows that
$\rho_1^{(i)}:\pi_1(X_1)\to\overline{\mathbb Q}_\ell^\ast$
corresponds to the sheaf $\mathscr K_{\chi^{q^i}, \sigma^{-i}(f)}$.
We also have
$$\sigma^{i}\Big(\frac{g\sigma^{-i}(\xi)}{\sigma^{-i}(\xi)}\Big)=
\Big(\frac{g\sigma^{-i}(\xi)}{\sigma^{-i}(\xi)}\Big)^{q^i}=\frac{g
\Big(\big(\sigma^{-i}(\xi)\big)^{q^i}\Big)}{\big(\sigma^{-i}(\xi)\big)^{q^i}},$$
and hence
$$\rho_1^{(i)}(g)=\chi^{-1}\left(\frac{g
\Big(\big(\sigma^{-i}(\xi)\big)^{q^i}\Big)}{\big(\sigma^{-i}(\xi)\big)^{q^i}}\right)$$
for any $g\in \pi_1(X_1)$. Note that
$\big(\sigma^{-i}(\xi)\big)^{q^i}$ is a root of the polynomial
$T^{q^d-1}-f(t_1^{q^i},\ldots, t_n^{q^i})$. This shows that
$\rho_1^{(i)}:\pi_1(X_1)\to\overline{\mathbb Q}_\ell^\ast$ also
corresponds to the sheaf $\mathscr K_{\chi, f(t_1^{q^i},\ldots,
t_n^{q^i})}$. By Proposition 1.2 (ii), $\pi^\ast(\otimes\hbox{-{\rm
Ind}}(\mathscr K_{\chi,f}))$ corresponds to the character
$\prod_{i=0}^{d-1} \rho_1^{(i)}$. So we get the following.

\begin{proposition}  Let $\chi:\mathbb F_{q^d}^\ast\to \overline{\mathbb
Q}_\ell^\ast$ be a multiplicative character, let $f\in\mathbb
F_{q^d}(t_1,\ldots, t_n)$, and let $X_0$ be an open subscheme of
$\mathbb A^n_{\mathbb F_q}$ so that its inverse image in $\mathbb
A^n_{\mathbb F_{q^d}}$ is contained in the complement of the union
of hypersurfaces $\{f=0\}\cup\{f=\infty\}$. Let $X_1$ be the base
change of $X_0$ from $\mathbb F_q$ to $\mathbb F_{q^d}$, and let
$\pi:X_1\to X_0$ be the projection. Then on $X_1$, we have
isomorphisms
$$\pi^\ast(\otimes\hbox{-{\rm Ind}}(\mathscr K_{\chi,f}))\cong
\bigotimes_{i=0}^{d-1} \mathscr K_{\chi^{q^i}, \sigma^{-i}(f)}\cong
\mathscr K_{\chi, \prod_{i=0}^{d-1}f(t_1^{q^i},\ldots,
t_n^{q^i})},$$ where $\sigma^{-i}(f)$ is the rational function
obtained from $f$ by taking the $q^i$-th root for each coefficient
of the numerator and the denominator of $f$.
\end{proposition}

\section{Tensor Induction of Artin-Schreier type sheaves}

The Artin-Schreier covering $$\wp:{\mathbb A}_{\mathbb F_q}^1\to
{\mathbb A}_{\mathbb F_q}^1,\quad x\mapsto x^q-x$$ defines an
$\mathbb F_q$-torsor
$$0\to \mathbb F_q \to {\mathbb A}_{\mathbb F_q}^1
\stackrel {\wp}\to {\mathbb A}_{\mathbb F_q}^1\to 0.$$ Let $\psi:
\mathbb F_q\to \overline {\mathbb Q}_\ell^\ast$ be a nontrivial
additive character. Pushing-forward this torsor by $\psi^{-1}$, we
get a lisse $\overline{\mathbb Q}_\ell$-sheaf ${\mathscr L}_\psi$ of
rank $1$ on ${\mathbb A}_{\mathbb F_q}^1$, which we call the
\emph{Artin-Schreier sheaf}. For any $x\in\mathbb A_{\mathbb
F_q}^1(\mathbb F_{q^m})=\mathbb F_{q^m}$, we have
$$\mathrm{Tr}(F_x,\mathscr L_{\psi,\bar x})=\psi(\mathrm{Tr}_{\mathbb F_{q^m}/\mathbb
F_q}(x)),$$ where $\mathrm{Tr}_{\mathbb F_{q^m}/\mathbb F_q}$
denotes the trace for the field extension $\mathbb F_{q^m}/\mathbb
F_q$.

Let $f\in\Gamma(X_0,\mathcal O_{X_0})$ be a section of the structure
sheaf $\mathcal O_{X_0}$. Then the $\mathbb F_q$-algebra
homomorphism
$$\mathbb F_q[t]\to \Gamma(X_0,\mathcal O_{X_0}),\quad t\mapsto f$$
defines an $\mathbb F_q$-morphism $X_0\to\mathbb A_{\mathbb F_q}^1$
which we still denote by $f$. Denote $f^\ast \mathscr L_\psi$ by
$\mathscr L_{\psi,f}$. It is a lisse $\overline{\mathbb
Q}_\ell$-sheaf of rank one on $X_0$. For any $x\in X_0(\mathbb
F_{q^m})$, we have
$$\mathrm{Tr}(F_x,(\mathscr L_{\chi,f})_{\bar x})=\psi\big(\mathrm{Tr}_{\mathbb F_{q^m}/\mathbb
F_q}(f(x))\big).$$ For any $f_1,f_2\in\Gamma(X_0,\mathcal
O_{X_0}^\ast)$, we have
$$\mathscr L_{\psi,f_1+f_2}\cong \mathscr L_{\psi,f_1}\otimes
\mathscr L_{\psi,f_2}.$$ Confer \cite[Sommes trig. 1.4-8]{SGA41/2}.

Suppose that $X_0$ is a normal scheme. Then $\mathscr L_{\psi,f}$
can be described as follows. Let $K$ be the function field of $X_0$.
Fix a separable closure of $\overline {K}$ of $K(X_0)$, and let the
base point of $X_0$ be the canonical morphism
$\mathrm{Spec}\,\overline K\to X_0$. Let $\zeta\in \overline K$ be a
root of the polynomial $T^q-T-f.$ Then any root of $T^q-T-f$ is of
the form $\zeta+a$ for some $a\in\mathbb F_q$, $K[\zeta]$ is a
Galois extension and the normalization of $X_0$ in $K[\zeta]$ is
etale over $X_0$. It follows that we have a canonical epimorphism
$$\pi_1(X_0)\to \mathrm{Gal}(K[\zeta]/K).$$ On the other hand, we have
a canonical monomorphism
$$\mathrm{Gal}(K[\zeta]/K)\hookrightarrow \mathbb F_q,\quad
g\mapsto g(\zeta)-\zeta.$$ The representation corresponding to
$\mathscr L_{\psi,f}$ is the composite
$$\pi_1(X_0)\to \mathrm{Gal}(K[\zeta]/K) \hookrightarrow \mathbb F_q
\stackrel{\psi^{-1}}\to \overline {\mathbb Q}_\ell^\ast.$$

Let $\mathbb F_{q^d}$ be an extension of $\mathbb F_q$ of degree
$d$, let $\psi:\mathbb F_{q^d}\to \overline{\mathbb Q}_\ell^\ast$ be
a nontrivial additive character, let $f\in\mathbb
F_{q^d}(t_1,\ldots, t_n)$ be a rational function, and let $X_0$ be
an open subscheme of $\mathbb A^n_{\mathbb F_q}$ so that its inverse
image is contained in the complement of the hypersurface $f=\infty$.
Let $X_1$ be the base change of $X_0$ from ${\mathbb F_q}$ to
$\mathbb F_{q^d}$, and let $\pi:X_1\to X_0$ be the projection.
Consider the sheaf $\mathscr L_{\psi,f}$ on $X_1$. Let us describe
the character $\pi_1(X_1)\to\overline{\mathbb Q}_\ell^\ast$
corresponding to the rank one lisse sheaf
$\pi^\ast(\otimes\hbox{-Ind}(\mathscr L_{\psi,f}))$ on $X_1$. The
function field of $X_0$ (resp. $X_1$) is the rational function field
$K_0=\mathbb F_q(t_1,\ldots, t_n)$ (resp. $K_1=\mathbb
F_{q^d}(t_1,\ldots, t_n)$). Let $\zeta$ be a root of the polynomial
$T^{q^d}-T-f$ in a separable closure $\overline K_1$ of
$K_1={\mathbb F_{q^d}(t_1,\ldots, t_n)}$. Then $\mathscr L_{\psi,f}$
corresponds to the character
$$\pi_1(X_1)\to
\mathrm{Gal}(K_1[\zeta]/K_1)\hookrightarrow \mathbb
F_{q^d}\stackrel{\psi^{-1}}\to \overline{\mathbb Q}_\ell^\ast.$$
Denote this character by $\rho_1:\pi_1(X_1)\to\overline{\mathbb
Q}_\ell^\ast$. More explicitly, we have
$$\rho_1(g)=\psi^{-1}\big(g(\zeta)-\zeta\big)$$ for any $g\in
\pi_1(X_1)$. We have canonical isomorphisms
$$\mathrm{Aut}(X_1/X_0)^\circ\cong\mathrm{Gal}\Big(\mathbb F_{q^d}(t_1,\ldots,
t_n)/\mathbb F_{q}(t_1,\ldots, t_n)\Big)\cong \mathrm{Gal}(\mathbb
F_{q^d}/\mathbb F_q).$$ Choose $\sigma\in \pi_1(X_0)$ so that its
image in $\mathrm{Gal}\Big(\mathbb F_{q^d}(t_1,\ldots, t_n)/\mathbb
F_{q}(t_1,\ldots, t_n)\Big)$ is given by
$$\sigma(a)=a^q \hbox { for any } a\in\mathbb F_{q^d},\quad
\sigma(t_i)=t_i\; (i=1, \ldots, n).$$ For any $0\leq i\leq d-1$, let
$\rho_1^{(i)}:\pi_1(X_1)\to\overline{\mathbb Q}_\ell^\ast$ be the
character defined by
$$\rho_1^{(i)}(g)=\rho_1(\sigma^ig\sigma^{-i})$$ for any $g\in
\pi_1(X_1)$. We have
\begin{eqnarray*}
\rho_1^{(i)}(g)&=&\rho_1(\sigma^{i}g\sigma^{-i})\\
&=& \psi^{-1}\big(\sigma^{i}g\sigma^{-i}(\zeta)-\zeta\big)\\
&=&\psi^{-1}\sigma^{i}\big(g\sigma^{-i}(\zeta)-\sigma^{-i}(\zeta)\big).
\end{eqnarray*}
Note that $\sigma^{-i}(\zeta)$ and $g\sigma^{-i}(\zeta)$ are roots
of the polynomial $T^{q^d}-T-\sigma^{-i}(f)$, and hence
$g\sigma^{-i}(\zeta)-\sigma^{-i}(\zeta)$ lies in $\mathbb F_{q^d}$.
We thus have
$$\sigma^{i}\big(g\sigma^{-i}(\zeta)-\sigma^{-i}(\zeta)\big)=
\big(g\sigma^{-i}(\zeta)-\sigma^{-i}(\zeta)\big)^{q^i}.$$ So we have
$$\rho_1^{(i)}(g)=\psi^{-1}\Big(\big(g\sigma^{-i}(\zeta)-\sigma^{-i}(\zeta)\big)^{q^i}\Big)$$
for any $g\in \pi_1(X_1)$. This shows that
$\rho_1^{(i)}:\pi_1(X_1)\to\overline{\mathbb Q}_\ell^\ast$
corresponds to the sheaf $\mathscr L_{\psi\circ \sigma^i,
\sigma^{-i}(f)}$. By Proposition 1.2 (ii),
$\pi^\ast(\otimes\hbox{-{\rm Ind}}(\mathscr L_{\psi,f}))$
corresponds to the character $\prod_{i=0}^{d-1} \rho_1^{(i)}$. So we
get the following.

\begin{proposition}  Let $\psi:\mathbb F_{q^d}\to \overline{\mathbb
Q}_\ell^\ast$ be a nontrivial additive character, let $f\in\mathbb
F_{q^d}(t_1,\ldots, t_n)$ be a rational function, and let $X_0$ be
an open subscheme of $\mathbb A^n_{\mathbb F_q}$ so that its inverse
image in $\mathbb A^n_{\mathbb F_{q^d}}$ is contained in the
complement of the hypersurface $f=\infty$. Let $X_1$ be the base
change of $X_0$ from $\mathbb F_q$ to $\mathbb F_{q^d}$, and let
$\pi:X_1\to X_0$ be the projection. Then on $X_1$, we have an
isomorphism
$$\pi^\ast(\otimes\hbox{-{\rm Ind}}(\mathscr L_{\psi,f}))\cong
\bigotimes_{i=0}^{d-1} \mathscr L_{\psi\circ \sigma^i,
\sigma^{-i}(f)},$$ where $\sigma^{-i}(f)$ is the rational function
obtained from $f$ by taking the $q^i$-th root for each coefficient
of the numerator and the denominator of $f$, and $\psi\circ\sigma^i$
is the additive character $a\mapsto \psi(a^{q^i})$ for any
$a\in\mathbb F_{q^d}$.
\end{proposition}

\begin{remark} Let $\psi_q:\mathbb
F_q\to\overline{\mathbb Q}_\ell^\ast$ be a nontrivial additive
character of $\mathbb F_q$, and let $\psi:\mathbb
F_{q^d}\to\overline{\mathbb Q}_\ell^\ast$ be the additive character
$\psi=\psi_q\circ\mathrm{Tr}_{\mathbb F_{q^d}/\mathbb F_q}$. In the
notation of Proposition 3.1, we have $\psi\circ \sigma^i=\psi$ and
\begin{eqnarray*}
\pi^\ast(\otimes\hbox{-{\rm Ind}}(\mathscr L_{\psi,f}))&\cong&
\bigotimes_{i=0}^{d-1}  \mathscr L_{\psi\circ \sigma^i,
\sigma^{-i}(f)}\\
&\cong& \bigotimes_{i=0}^{d-1}  \mathscr L_{\psi, \sigma^{-i}(f)}\\
&\cong& \mathscr L_{\psi, \mathrm{Tr}_{\mathbb F_{q^d}/\mathbb
F_q}(f)},
\end{eqnarray*}
where $\mathrm{Tr}_{\mathbb F_{q^d}/\mathbb
F_q}(f)=\sum_{i=0}^{d-1}\sigma^{i}(f)$.
\end{remark}

\section{Character Sums on Curves}

In this section, $X_0$ is a smooth geometrically connected curve
over $\mathbb F_q$ of genus $g$, $\overline X_0$ is the smooth
compactification of $X_0$, $X$ and $\overline X$ are the base
changes of $X_0$ and $\overline X_0$ from $\mathbb F_q$ to $\mathbb
F$, respectively. Let us recall the general method of studying
character sums on $X_0$ using $\ell$-adic cohomology theory. Let
$\mathscr L_0$ be a lisse $\overline {\mathbb Q}_\ell$-sheaf on
$X_0$ of rank $r$. Suppose the inverse image $\mathscr L$ on $X$ has
neither nonzero constant subsheaf nor nonzero constant quotient
sheaf. Then by \cite[1.4.1]{D1}, we have
$$H^0_c(X,\mathscr L)=0,\quad
H^2_c(X,\mathscr L)=0.$$ By the Grothendieck-Ogg-Shafarevich formula
(\cite[X 7.1]{SGA5}), we have
\begin{eqnarray*}
\mathrm{dim}\, H^1_c(X, \mathscr L)&=&-\sum_{i=0}^2 (-1)^i
\mathrm{dim}\, H^i_c(X, \mathscr L)\\
&=&r \big(2g-2+\#(\overline X-X)\big)+\sum_{x\in \overline
X-X}\mathrm{sw}_x(\mathscr L),
\end{eqnarray*}
where $\mathrm{sw}_x$ denotes the Swan conductor at $x$ for any
Zariski closed point $x$ in $\overline X$. Suppose furthermore that
$\mathscr L_0$ is punctually pure of weight $w$. (Confer
\cite[D\'efinition 1.2.2]{D1}). Then by \cite[Corollarie 3.3.4]{D1},
all the eigenvalues of $F$ on $H^1_c(X, \mathscr L)$ have
Archimedean absolute value $\leq q^{\frac{w+1}{2}}$. By the
Grothendieck trace formula \cite[Rapport 3.2]{SGA41/2}, we have
\begin{eqnarray*}
S_m(X_0,\mathscr L_0)&=&\sum_{x\in X_0(\mathbb
F_{q^m})}\mathrm{Tr}(F_x,\mathscr L_{0,\bar
x})\\
&=&\sum_{i=0}^2 (-1)^i\mathrm{Tr}\big(F^m, H^i_c(X,
\mathscr L)\big)\\
&=&-\mathrm{Tr}\big(F^m, H^1_c(X, \mathscr L)\big).
\end{eqnarray*} So we have
$$|S_m(X_0,\mathscr L_0)|\leq \Big(r\big(2g-2+\#(\overline X-X)\big)+
\sum_{x\in \overline X-X}\mathrm{sw}_x(\mathscr
L)\Big)q^{\frac{m(w+1)}{2}}.$$ Finally since $\mathscr L_0$ is
punctually pure, $\mathscr L$ is semi-simple (\cite[Th\'eor\`eme
3.4.1 (iii)]{D1}). The condition that $\mathscr L$ has  has neither
nonzero constant subsheaf nor nonzero constant quotient sheaf is
equivalent to the condition that $\mathscr L$ has no nonzero
constant subsheaf. We thus get the following.

\begin{proposition} Let $X_0$ be a smooth geometrically connected curve
over $\mathbb F_q$ of genus $g$, $\overline X_0$ its smooth
compactification, and $\mathscr L_0$ a lisse $\overline {\mathbb
Q}_\ell$-sheaf on $X_0$ of rank $r$. Suppose $\mathscr L_0$ is
punctually pure of weight $w$, and its inverse image $\mathscr L$ on
$X$ has no nonzero constant subsheaf. Then we have
\begin{eqnarray*}
| S_m(X_0,\mathscr L_0)|&=&\bigg|\sum_{x\in X_0(\mathbb
F_{q^m})}\mathrm{Tr}(F_x,\mathscr
L_{0, \bar x})\bigg|\\
&\leq&\Big(r\big(2g-2+\#(\overline X-X)\big)+ \sum_{x\in \overline
X-X}\mathrm{sw}_x(\mathscr L)\Big)q^{\frac{m(w+1)}{2}}.
\end{eqnarray*}
\end{proposition}

\begin{theorem} Let $X_0$ be a smooth geometrically connected curve
over $\mathbb F_q$ of genus $g$, $\overline X_0$ its smooth
compactification, $\mathbb F_{q^d}$ an extension of $\mathbb F_q$ of
degree $d$, $X_1=X_0\otimes_{\mathbb F_q}\mathbb F_{q^d}$, and
$\mathscr L_1$ a lisse $\overline {\mathbb Q}_\ell$-sheaf on $X_1$
of rank $r$. Suppose $\mathscr L_1$ is punctually pure of weight
$w$, and suppose the inverse image of $\otimes\hbox{-{\rm
Ind}}(\mathscr L_1)$ on $X$ has no nonzero constant subsheaf. Let
$$S_1'(X_1,\mathscr L_1)=\sum_{x\in X_0(\mathbb
F_q)}\mathrm{Tr}(F_{x_1},\mathscr L_{1,\bar x_1}),$$ where
$x_1:\mathrm{Spec}\,\mathbb F_{q^d}\to X_1$ is the base change of
$x:\mathrm{Spec}\,\mathbb F_q\to X_0$ for any $x\in X_0(\mathbb
F_q)$. Then we have
$$
| S_1'(X_1,\mathscr L_1)|\leq \Big(\big(2g-2+\#(\overline X-X)\big)+
d\sum_{x\in \overline X-X}\mathrm{sw}_x(\mathscr L)\Big)\cdot r^d
\cdot q^{\frac{dw+1}{2}}.$$
\end{theorem}

\begin{proof} By Proposition 1.2 (i) and Grothendieck's trace formula (\cite[Rapport 3.2]{SGA41/2}), we have
\begin{eqnarray*}
S_1'(X_1,\mathscr L_1)&=&\sum_{x\in X_0(\mathbb
F_q)}\mathrm{Tr}\big(F_x,(\otimes\hbox{-{\rm Ind}}(\mathscr
L_1))_{\bar
x}\big)\\
&=&\sum_{i=0}^2(-1)^i \mathrm{Tr}\big(F,H_c^i(X,\otimes\hbox{-{\rm
Ind}}(\mathscr L_1))\big).
\end{eqnarray*}
Since the inverse image of $\otimes\hbox{-{\rm Ind}}(\mathscr L_1)$
on $X$ has neither nonzero constant subsheaf nor nonzero constant
quotient sheaf, we have $H_c^i(X,\otimes\hbox{-{\rm Ind}}(\mathscr
L_1))=0$ for $i\not=1$. So we have $$S_1'(X_1,\mathscr
L_1)=-\mathrm{Tr}\big(F,H_c^1(X,\otimes\hbox{-{\rm Ind}}(\mathscr
L_1))\big).$$ One can verify $\otimes\hbox{-{\rm Ind}}(\mathscr
L_1)$ is punctually pure of weight $dw$. By \cite[Corollarie
3.3.4]{D1}, all eigenvalues of $F$ on $H_c^1(X,\otimes\hbox{-{\rm
Ind}}(\mathscr L_1))$ have Archimedean absolute value $\leq
q^{\frac{dw+1}{2}}$. So we have
$$|S_1'(X_1,\mathscr
L_1)|\leq \big(\mathrm{dim}\,H_c^1(X,\otimes\hbox{-{\rm
Ind}}(\mathscr L_1))\big)q^{\frac{dw+1}{2}}.$$ On the other hand, by
the Grothendieck-Ogg-Shafarevich formula (\cite[X 7.1]{SGA5}), we
have
\begin{eqnarray*}
\mathrm{dim}\,H_c^1(X,\otimes\hbox{-{\rm Ind}}(\mathscr
L_1))&=&-\sum_{i=0}^2
(-1)^i\mathrm{dim}\,H_c^i(X,\otimes\hbox{-{\rm Ind}}(\mathscr L_1))\\
&=& r^d\big(2g-2+\#(\overline X-X)\big)+ \sum_{x\in \overline
X-X}\mathrm{sw}_x(\otimes\hbox{-{\rm Ind}}(\mathscr L_1)).
\end{eqnarray*}
To prove our theorem, it suffices to show that
$$\sum_{x\in \overline
X-X}\mathrm{sw}_x(\otimes\hbox{-{\rm Ind}}(\mathscr L_1))\leq dr^d
\sum_{x\in \overline X-X}\mathrm{sw}_x(\mathscr L_1).$$ Choose
$\sigma\in\pi_1(X_0)$ so that its image under the canonical
epimorphism $\pi_1(X_0)\to\mathrm{Gal}(\mathbb F_{q^d}/\mathbb F_q)$
is the Frobenius substitution. Let
$\rho_1:\pi_1(X_1)\to\mathrm{GL}(r, \overline{\mathbb Q}_\ell)$ be
the representation corresponding to $\mathscr L_1$ and let
$\rho_1^{(i)}:\pi_1(X_1)\to\mathrm{GL}(r,\overline{\mathbb Q}_\ell)$
$(0\leq i\leq d-1)$ be the representations defined by
$$\rho_1^{(i)}(g)=\rho_1(\sigma^ig\sigma^{-i})$$ for any
$g\in\pi_1(X_1)$. Then by Proposition 1.2 (ii),
$\pi^\ast(\otimes\hbox{-{\rm Ind}}(\mathscr L_1))$ corresponds to
the representation $\otimes_{i=0}^{d-1}\rho_1^{(i)}$, where
$\pi:X_1\to X_0$ is the projection. Note that $\sigma$ induce an
automorphism of $X$ over $X_0$, and an automorphism of $\overline X$
over $\overline X_0$. For any $x\in \overline X-X$ and any $0\leq
i\leq d-1$, we have
$$\mathrm{sw}_x(\rho_1^{(i)})=\mathrm{sw}_{\sigma^i(x)}(\rho_1).$$
Combined with Lemma 4.3 below, we get
\begin{eqnarray*}
\mathrm{sw}_x(\otimes\hbox{-{\rm Ind}}(\mathscr L_1))&=&
\mathrm{sw}_x\Big(\otimes_{i=0}^{d-1}\rho_1^{(i)}\Big)\\
&\leq& r^d \max\{\mathrm{sw}_x(\rho_1^{(i)})|0\leq i\leq d-1\}\\
&=& r^d \max\{\mathrm{sw}_{\sigma^i(x)}(\rho_1)|0\leq i\leq d-1\}\\
&=& r^d \max\{\mathrm{sw}_{\sigma^i(x)}(\mathscr L_1)|0\leq i\leq
d-1\}\\
&\leq&r^d \sum_{i=0}^{d-1} \mathrm{sw}_{\sigma^i(x)}(\mathscr L_1).
\end{eqnarray*}
Hence
\begin{eqnarray*}
\sum_{x\in \overline X-X}\mathrm{sw}_x(\otimes\hbox{-{\rm
Ind}}(\mathscr L_1))&\leq& r^d \sum_{x\in \overline X-X}
\sum_{i=0}^{d-1}
\mathrm{sw}_{\sigma^i(x)}(\mathscr L_1)\\
&=&d r^d \sum_{x\in \overline X-X} \mathrm{sw}_x(\mathscr L_1).
\end{eqnarray*}
\end{proof}

\begin{lemma} Let $K$ be a local field whose residue field is of characteristic prime to $\ell$, and let
$U$ and $V$ be $\overline {\mathbb Q}_\ell$-representations of
$\mathrm{Gal}(\overline K/K)$. Then we have the following estimation
for the Swan conductor of tensor product
$$\mathrm{sw}(U\otimes V)\leq \mathrm{dim}(U\otimes V)\cdot \max(\mathrm
{sw}(U),\mathrm{sw}(V)).$$
\end{lemma}

\begin{proof} Let $G^{(\lambda)}$ $(\lambda\geq 0)$ be the filtration of
$\mathrm{Gal}(\overline K/K)$ in upper numbering, let $G^{(\lambda
+)}$ be the closure of $\bigcup_{\epsilon
>0}G^{(\lambda +\epsilon)}$, and let
$$U=\bigoplus_\lambda U_\lambda,\quad  V=\bigoplus_\lambda
V_\lambda$$ be the decompositions of $U$ and $V$ as representations
of the wild inertia subgroup $G^{(0+)}$ such that
\begin{eqnarray*}
U_\lambda^{G^{(\lambda+)}}=U_\lambda, &&
U_\lambda^{G^{(\lambda)}}=0\hbox { (if } \lambda>0),\\
V_\lambda^{G^{(\lambda+)}}=V_\lambda, &&
V_\lambda^{G^{(\lambda)}}=0\hbox { (if } \lambda>0).
\end{eqnarray*}
For any $\lambda, \mu\geq 0$, we have
$$(U_\lambda\otimes V_\mu)^{G^{(\max(\lambda,\mu)+)}}=U_\lambda\otimes
V_\mu.$$ On the other hand, we have
$$U\otimes V=\bigoplus_{\lambda,\mu } U_\lambda\otimes V_\mu.$$ It
follows that
\begin{eqnarray*}
\mathrm{sw}(U\otimes V)&\leq&
\sum_{\lambda,\mu}\max(\lambda,\mu)\mathrm{dim}(U_\lambda)\mathrm{dim}(V_\mu)\\
&\leq& \max(\mathrm {sw}(U),\mathrm{sw}(V)) \sum_{\lambda,\mu }
\mathrm{dim}(U_\lambda)\mathrm{dim}(V_\mu)\\
&=&\mathrm{dim}(U\otimes V)\cdot \max(\mathrm
{sw}(U),\mathrm{sw}(V)).
\end{eqnarray*}
\end{proof}

\begin{remark} In the proof of Theorem 4.2, suppose there exists closed
points $\infty_1,\ldots, \infty_m\in \overline X-X$ such that
$\sigma(\infty_j)=\infty_j$ $(i=1,\ldots, m)$. Then we have
\begin{eqnarray*}
\mathrm{sw}_{\infty_j}(\otimes\hbox{-{\rm Ind}}(\mathscr L_1))&=&
\mathrm{sw}_{\infty_j}\Big(\otimes_{i=0}^{d-1}\rho_1^{(i)}\Big)\\
&\leq& r^d \max\{\mathrm{sw}_{\infty_j}(\rho_1^{(i)})|0\leq i\leq d-1\}\\
&=& r^d \max\{\mathrm{sw}_{\sigma^i({\infty_j})}(\rho_1)|0\leq i\leq d-1\}\\
&=& r^d \mathrm{sw}_{{\infty_j}}(\mathscr L_1).
\end{eqnarray*}
The estimate in Theorem 4.2 can be improved as $$ |
S_1'(X_1,\mathscr L_1)|\leq\Big(\big(2g-2+\#(\overline X-X)\big)+
d\sum_{x\in \overline X-X,\; x\not=\infty_1,\ldots,
\infty_m}\mathrm{sw}_x(\mathscr L_1)+ \sum_{j=1}^m
\mathrm{sw}_{{\infty_j}}(\mathscr L_1)\Big)\cdot r^d \cdot
q^{\frac{dw+1}{2}}.$$
\end{remark}

From now on, let $\sigma$ be the automorphism of $\mathbb
F_{q^d}(t)$ defined by $\sigma(t)=t$ and $\sigma(a)=a^q$ for any
$a\in \mathbb F_{q^d}$. We will apply Theorem 4.2 to Kummer type and
Artin-Schreier type sheaves. Note that these sheaves have rank $1$.

\begin{theorem} Let $f(t)\in\mathbb F_{q^d}(t)$ be a rational function. Write
$f(t)=\prod_{j=1}^k f_j(t)^{n_j}$, where $f_j(t)\in\mathbb
F_{q^d}[t]$ are irreducible polynomials and $n_j$ are nonzero
integers. Let $\chi:\mathbb F_{q^d}^\ast\to\overline{\mathbb
Q}_\ell^\ast$ be a multiplicative character for ${\mathbb F}_{q^d}$.
Suppose that the rational function $\prod_{i=0}^{d-1}f(t^{q^i})$ is
not of the form $h(t)^{\mathrm{ord}(\chi)}$ in $\mathbb F(t)$, where
$\mathrm{ord}(\chi)$ is the smallest integer $d$ such that
$\chi^d=1$. Then we have
$$\bigg|\sum_{a\in {\mathbb F}_q,\; f(a)\not=0,\infty}
\chi(f(a))\bigg|  \leq \Big(d\sum_{j=1}^k \mathrm{deg}(f_j)-
1\Big)\sqrt{q}.
$$
\end{theorem}

\begin{proof} Let $X_0$ be the
complement in $\mathbb A^1_{\mathbb F_q}$ of the hypersurface
$\prod_{i=0}^{d-1}\prod_{j=1}^k\sigma^{-i}(f_j)=0$. Consider the
$\overline {\mathbb Q}_\ell$-sheaf $\mathscr K_{\chi, f}$ on
$X_1=X_0\otimes_{\mathbb F_q}\mathbb F_{q^d}$. We have
$$\sum_{a\in {\mathbb F}_q,\; f(a)\not=0,\infty}\chi(f(a))=\sum_{x\in X_0(\mathbb
F_q)}\mathrm{Tr}(F_{x_1},(\mathscr K_{\chi, f})_{\bar x_1}).$$ Note
that $\mathscr K_{\chi, f}$ is tamely ramified everywhere on
$\overline X=\mathbb P^1_{\mathbb F}$, and hence
$\mathrm{sw}_x(\mathscr K_{\chi, f})=0$ for all $x\in \overline
X-X=\mathbb P^1_{\mathbb F} -X_0\otimes_{\mathbb F_q}\mathbb F$. The
set $\overline X-X$ consists of $\infty$ and the roots of
$\prod_{i=0}^{d-1}\prod_{j=1}^k\sigma^{-i}(f_j)$. There are
$d\sum_{j=1}^k \mathrm{deg}(f_j)$ such roots. The genus of $X_0$ is
$0$, and $\mathscr K_{\chi,f}$ is punctually pure of weight $0$ and
of rank $1$. Our estimate follows directly from the estimate in
Theorem 4.2, provided that we can prove $\otimes\hbox{-{\rm
Ind}}(\mathscr K_{\chi, f})$ is geometrically nonconstant. Indeed,
by Proposition 2.1, we have $\pi^\ast(\otimes\hbox{-{\rm
Ind}}(\mathscr K_{\chi, f}))\cong \mathscr K_{\chi,
\prod_{i=0}^{d-1}f(t^{q^i})}$, where $\pi:X_1\to X_0$ is the
projection. Since $\prod_{i=0}^{d-1}f(t^{q^i})$ is not of the form
$h(t)^{\mathrm{ord}(\chi)}$ in $\mathbb F(t)$, the sheaf $\mathscr
K_{\chi, \prod_{i=0}^{d-1}f(t^{q^i})}$ is geometrically nonconstant.
\end{proof}

\begin{theorem} Let $f(t), g(t)\in\mathbb F_{q^d}(t)$ be rational functions. Write
$f(t)=\prod_{j=1}^k f_j(t)^{n_j}$, where $f_j(t)\in\mathbb
F_{q^d}[t]$ are irreducible polynomials and $n_j$ are nonzero
integers. Let $D_1=\sum_{j=1}^k\mathrm{deg}(f_j)$, let
$D_2=\max(\mathrm {deg}(g), 0)$, let $D_3$ be the degree of the
denominator of $g(t)$, and let $D_4$ be the sum of degrees of those
irreducible polynomials dividing the denominator of $g$ but distinct
from $f_j(t)$ $(j=1,\ldots, k)$. Let $\chi:\mathbb
F_{q^d}^\ast\to\overline{\mathbb Q}_\ell^\ast$ be a multiplicative
character of $\mathbb F_{q^d}$, and let $\psi_q:\mathbb
F_q\to\overline{\mathbb Q}_\ell^\ast$ be a nontrivial additive
character of $\mathbb F_q$. Suppose $\mathrm{Tr}_{\mathbb
F_{q^d}/\mathbb F_q}(g)=\sum_{i=0}^{d-1} \sigma^i(g)$ is not of the
form $r(t)^p-r(t)$ in ${\mathbb F}(t)$. Then we have the estimate
$$\bigg|
\sum_{a\in {\mathbb F}_q, \;f(a)\not= 0,\infty,\;
g(a)\not=\infty}\chi(f(a))\psi_q\Big(\mathrm{Tr}_{\mathbb
F_{q^d}/\mathbb F_q}(g(a))\Big)\bigg| \leq
(d(D_1+D_3+D_4)+D_2-1)\sqrt{q}.$$
\end{theorem}

\begin{proof} Let $u(t)$ be the denominator of $g(t)$,
let $X_0$ be the complement in $\mathbb A^1_{\mathbb F_q}$ of the
hypersurface $\prod_{i=0}^{d-1}\Big(\sigma^{-i}(u)
\prod_{j=1}^k\sigma^{-i}(f_j)\Big)=0$, and let $\mathscr L_1$ be the
$\overline {\mathbb Q}_\ell$-sheaf $\mathscr L_1=\mathscr K_{\chi,
f}\bigotimes \mathscr L_{\psi, g}$ on $X_1=X_0\otimes_{\mathbb
F_q}\mathbb F_{q^d}$, where $\psi=\psi_q\circ\mathrm{Tr}_{\mathbb
F_{q^d}/\mathbb F_q}$. Then we have
$$\sum_{a\in {\mathbb F}_q, \;f(a)\not= 0,\infty,\;
g(a)\not=\infty}\chi(f(a))\psi_q\Big(\mathrm{Tr}_{\mathbb
F_{q^d}/\mathbb F_q}(g(a))\Big)=\sum_{x\in X_0(\mathbb
F_q)}\mathrm{Tr}(F_{x_1},\mathscr L_{1,\bar x_1}).$$ Note that
$\mathscr K_{\chi, f}$ is tamely ramified everywhere on $\overline
X=\mathbb P^1_{\mathbb F}$, and has no effect on the Swan conductor
of $\mathscr L_1$. So for any $x\in \overline X-X=\mathbb
P^1_{\mathbb F}-X_0\otimes_{\mathbb F_q}\mathbb F$, we have
$$
\mathrm{sw}_x(\mathscr L_1)=\mathrm{sw}_x(\mathscr L_{\psi, g}).$$
Let $v_x$ be the valuation of the field $\mathbb F(t)$ corresponding
to the point $x\in\overline X-X$. If $x$ is not a pole of $g$, then
$\mathrm{sw}_x(\mathscr L_{\psi, g})=0$. Otherwise, we have
$$\mathrm{sw}_x(\mathscr L_{\psi, g})\leq -v_x(g).$$ It
follows that
$$\mathrm{sw}_\infty(\mathscr L_1)\leq D_2,\quad \sum_{x\in \overline
X-X,\; x\not=\infty}\mathrm{sw}_x(\mathscr L_1)\leq D_3.$$  The set
$\overline X-X$ consists of $\infty$ and roots of
$\prod_{i=0}^{d-1}\Big(\sigma^{-i}(u)
\prod_{j=1}^k\sigma^{-i}(f_j)\Big)$. It follows that
$$\#(\overline X-X)\leq
1+d(D_1+D_4).$$ The genus of $X_0$ is $0$, $\mathscr L_1$ is
punctually pure of weight $0$ and of rank $1$. Our estimate follows
directly from the estimate in Theorem 4.2 and Remark 4.4, provided
that we can prove $\otimes\hbox{-{\rm Ind}}(\mathscr L_1)$ is
geometrically nonconstant. Let $\pi:X_1\to X_0$ be the projection.
It suffices to prove $\pi^\ast(\otimes\hbox{-{\rm Ind}}(\mathscr
L_1))$ is geometrically nonconstant. By Propositions 1.2 (iii), 2.1,
3.1 and Remarks 3.2, we have
$$\pi^\ast(\otimes\hbox{-{\rm Ind}}(\mathscr L_1))\cong\Big(
\bigotimes_{i=0}^{d-1} \mathscr K_{\chi^{q^i},
\sigma^{-i}(f)}\Big)\bigotimes \mathscr L_{\psi,
\mathrm{Tr}_{\mathbb F_{q^d}/\mathbb F_q}(g)}.$$ Since
$\mathrm{Tr}_{\mathbb F_{q^d}/\mathbb F_q}(g)$ is not of the form
$r(t)^p-r(t)$ in ${\mathbb F}(t)$, the sheaf $\mathscr L_{\psi,
\mathrm{Tr}_{\mathbb F_{q^d}/\mathbb F_q}(g)}$ is geometrically
nonconstant, and it has wild ramification at poles of
$\mathrm{Tr}_{\mathbb F_{q^d}/\mathbb F_q}(g)$, whereas
$\bigotimes_{i=0}^{d-1} \mathscr K_{\chi^{q^i}, \sigma^{-i}(f)}$ is
at worst tamely ramified. So $(\bigotimes_{i=0}^{d-1}\mathscr
K_{\chi^{q^i}, \sigma^{-i}(f)})\bigotimes \mathscr L_{\psi,
\mathrm{Tr}_{\mathbb F_{q^d}/\mathbb F_q}(g)}$ is geometrically
nonconstant.
\end{proof}

\begin{remark} Let $f_1(t),\ldots, f_m(t),g_1(t),\ldots, g_n(t)\in\mathbb F_{q^d}(t)$ be rational functions,
let $\chi_1,\ldots, \chi_m:\mathbb F_{q^d}^\ast\to\overline{\mathbb
Q}_\ell^\ast$ be multiplicative characters, and let $\psi_1,\ldots,
\psi_n:\mathbb F_{q^d}\to\overline{\mathbb Q}_\ell^\ast$ be additive
characters. Then the sum
$$\sum_{a\in {\mathbb F}_q, \;f_i(a)\not= 0,\infty,\;
g_j(a)\not=\infty}\chi_1(f_1(a))\cdots
\chi_m(f_m(a))\psi_1(g_1(a))\cdots \psi_n(g_n(a))$$ can be reduced
to the form of the sum in Theorem 4.6. Indeed, let $\chi:\mathbb
F_{q^d}^\ast\to\overline{\mathbb Q}_\ell^\ast$ be a multiplicative
character of order $q^d-1$. We can find integers $k_i$ $(i=1,\ldots,
m)$ such that $\chi_i=\chi^{k_i}$. Then we have $\chi_i(f_i(a))=
\chi(f_i(a)^{k_i})$. Let $f(t)=\prod_{i=1}^m f_i(t)^{k_i}$. We can
find $\lambda_j\in\mathbb F_{q^d}$ $(j=1,\ldots, n)$ such that
$\psi_j(a)=\psi_q(\mathrm{Tr}_{\mathbb F_{q^d}/\mathbb
F_q}(\lambda_j a))$ for any $a\in\mathbb F_{q^d}$. Let
$g(t)=\sum_{j=1}^n \lambda_j g_j$. Then we have
\begin{eqnarray*}
&& \sum_{a\in {\mathbb F}_q, \;f_i(a)\not= 0,\infty,\;
g_j(a)\not=\infty}\chi_1(f_1(a))\cdots
\chi_m(f_m(a)))\psi_1(g_1(a))\cdots \psi_n(g_n(a))\\
&&=\sum_{a\in {\mathbb F}_q, \;f_i(a)\not= 0,\infty,\;
g_j(a)\not=\infty}\chi(f(a)) \psi_q\Big(\mathrm{Tr}_{\mathbb
F_{q^d}/\mathbb F_q}(g(a))\Big).
\end{eqnarray*}
\end{remark}

\end{document}